\documentclass[12pt]{article}
\usepackage{tikz}
\usetikzlibrary{calc}
\usepackage{bez123,calc,curves,ebezier,epic,eepic,graphicx,multiply,rotating}

\everymath{\displaystyle}
\usepackage{graphicx}
\usepackage{pst-all}
\usepackage{color}
\usepackage{amssymb}
\usepackage{amsmath}
\newtheorem{prethm}{{\bf Theorem}}

\newtheorem{prelemma}{{\bf Lemma}}

\newtheorem{preex}{{\bf Example}}

\newtheorem{preprop}{{\bf Proposition}}

\newenvironment{prop}{\begin{preprop}{\hspace{-0.5em}{\bf .}}}{\end{preprop}}
\newtheorem{precor}{{\bf Corollary}}

\newtheorem{preremark}{{\bf Remark}}

\newtheorem{preprob}{{\bf Problem}}

\newtheorem{predefin}{{\bf Definition}}

\newtheorem{preconj}{{\bf Conjecture}}

\newenvironment{conj}{\begin{preconj}{\hspace{-0.7
               em}{\bf.}}}{\end{preconj}}
\newtheorem{preprobb}{{\bf Problem}}

\newtheorem{prelem}{{\bf Theorem}}

\newenvironment{proof}{{\bf Proof.}\rm }{\hfill{$\Box$}}

\newtheorem{presolution}{{\bf Solution.}}

\def\newpic#1{}

\title{\vspace{-2.51cm}\Large\bf A note concerning the Grundy and ${\rm b}$-chromatic number of graphs}

\author{\large\bf Manouchehr Zaker\footnote{mzaker@iasbs.ac.ir}
	\vspace{5mm}\\
	Department of Mathematics,\\
	Institute for Advanced Studies in Basic Sciences,\\
	Zanjan 45137-66731, Iran\\
}

\date{}
\begin{document}
\maketitle
\begin{abstract}
\noindent The Grundy number of a graph $G$ is the maximum number of colors used by the First-Fit coloring of $G$ and is denoted by $\Gamma(G)$. Similarly, the ${\rm b}$-chromatic number ${\rm{b}}(G)$ of $G$ expresses the worst case behavior of another well-known coloring procedure i.e. color-dominating coloring of $G$. We obtain some families of graphs $\mathcal{F}$ for which there exists a function $f(x)$ such that $\Gamma(G)\leq f({\rm{b}}(G))$, for each graph $G$ from the family. Call any such family $(\Gamma,b)$-bounded family. We conjecture that the family of ${\rm b}$-monotone graphs is $(\Gamma,b)$-bounded and validate the conjecture for some families of graphs.
\end{abstract}




\section{Introduction}

\noindent This note deals only with undirected graphs without any loops or multiple edges. By a Grundy coloring of a graph $G$ we mean any partition of $V(G)$ into independent subsets $C_1, \ldots, C_k$ such that for each $i,j \in \{1, \ldots, k\}$ with $i<j$, each vertex in $C_j$ has a neighbor in $C_i$. The maximum such value $k$ is called the Grundy number (also called First-Fit chromatic number) and denoted by $\Gamma(G)$ (also by $\chi_{\sf FF}(G)$). It can be observed that $\Gamma(G)$ is equal to the maximum number of colors used by the First-Fit (greedy) coloring procedure in the graph $G$ \cite{Z1}. The Grundy number and First-Fit coloring of graphs are important research areas in chromatic and algorithmic graph theory with full of papers e.g. \cite{CS, GL, KP, KPT, TWHZ, Z1, Z2}.

\noindent By a color-dominating coloring of $G$ we mean any partition of $V(G)$ into independent subsets $C_1, \ldots, C_k$ such that for each $i$, the class $C_i$ contains a vertex say $v$ such that $v$ has a neighbor in any other class $C_j$, $j\not=i$. Denote by ${\rm{b}}(G)$ (also denoted by $\varphi(G)$) the maximum number of colors used in any color-dominating coloring of $G$. It can be easily seen that ${\rm{b}}(G)\leq \Delta(G)+1$ and under some conditions the equality holds, e.g. $d$-regular graphs with at least $2d^3$ vertices \cite{CJ}. An algorithmic interpretation of ${\rm{b}}(G)$ is that it expresses the worst case behavior of the following coloring procedure. In any proper coloring $C$ of a graph $G$, a vertex $v$ is said to be a color-dominating vertex if it has a neighbor with any other color except the color of $v$. Let $C$ be any arbitrary proper coloring of $G$ and $C_i$ be a color class in $C$.
If $C_i$ does not contain any color-dominating vertex then each vertex of $C_i$ can be removed from $C_i$ and transferred to another suitable class. By this technique the class $C_i$ is totally removed and number of colors is decreased by one. We repeat this method for all remaining color classes until we obtain a color-dominating coloring. Obviously, the final number of colors is at most ${\rm{b}}(G)$. The ${\rm b}$-chromatic number of graphs introduced in \cite{IM} and widely studied in the literature \cite{KZ, KTV}. For a recent survey on ${\rm b}$-chromatic number see \cite{JP}. A useful graph parameter relating to ${\rm b}$-chromatic number of a graph $G$ with non-increasing degree sequence $d_1\geq d_2 \geq \ldots \geq d_n$ is $m(G):=\max \{i: d_i\geq i-1\}$. It is known that ${\rm{b}}(G)\leq m(G)$ and for trees $T$, ${\rm{b}}(T)\geq m(T)-1$ \cite{IM}. In this paper, a graph $G$ is called ${\rm b}$-monotone if for each induced subgraph $H$ of $G$ we have ${\rm{b}}(H)\leq {\rm{b}}(G)$.

\noindent A first natural inquiry concerning the comparison of Grundy and ${\rm b}$-chromatic numbers is to explore and generate families of graphs $\{G_n\}_{n\geq 1}$ and $\{H_n\}_{n\geq 1}$ such that ${\rm{b}}(G_n)-\Gamma(G_n)\rightarrow \infty$ and $\Gamma(H_n)-{\rm{b}}(H_n)\rightarrow \infty$. Based on the results of this paper, both of the above-mentioned situations may happen in the universe of graphs. But the first situation (i.e. families with bounded Grundy number and unbounded ${\rm b}$-chromatic number) is more likely to happen because these families are more accessible.

\noindent The concept of $(\chi_{\sf FF}, \omega)$-boundedness was introduced by Gy{\'a}rf{\'a}s and Lehel in \cite{GL}. Denote the size of a maximum clique in $G$ by $\omega(G)$. A family $\mathcal{F}$ is called $(\chi_{\sf FF}, \omega)$-bounded if there exits a function $f(x)$ such that $\Gamma(G)\leq f(\omega(G))$ for each $G$ from the family. Some $(\chi_{\sf FF}, \omega)$-bounded families were obtained in \cite{GL,KP,KPT,Z2}. In the next section we introduce $(\Gamma, b)$-bounded families.

\section{$(\Gamma, b)$-bounded families of graphs}

\noindent We say a family $\mathcal{F}$ is $(\Gamma, b)$-bounded if
there exists a function $f(x)$ such that $\Gamma(G)\leq f({\rm{b}}(G))$, for each graph $G$ from the family. Note that any $(\chi_{\sf FF}, \omega)$-bounded family is also $(\Gamma, b)$-bounded. Also any family of graphs satisfying ${\rm{b}}(G)=\Delta(G)+1$ is $(\Gamma, b)$-bounded. Some of such families were obtained in \cite{KTV} and reported in \cite{JP}. With a similar manner we can define $(b, \Gamma)$-bounded families. We can easily obtain a sequence of trees $T_n$ such that $\Gamma(T_n)\leq 3$ for each $n$, but ${\rm{b}}(T_n)\rightarrow \infty$. In fact, we may consider $T_n$ as a path with sufficiently large length and sufficiently many leaves attached to the vertices of the path. In this note we concentrate on $(\Gamma, b)$-bounded families. The following proposition is useful.

\begin{prop}
A family $\mathcal{F}$ is $(\Gamma, b)$-bounded if and only if for any sequence $\{G_n\}_{n\geq 1}$ from $\mathcal{F}$, $\Gamma(G_n)\rightarrow \infty$ implies ${\rm{b}}(G_n)\rightarrow \infty$.
\end{prop}

\noindent \begin{proof}
If $\mathcal{F}$ is $(\Gamma, b)$-bounded then the assertion trivially holds. To prove the other side, note that any infinite family of graphs is countable, so write
$\mathcal{F}=\{G_n\}_{n\geq 1}$. If necessary use a relabeling and assume that $\{\Gamma(G_n)\}_{n\geq 1}$ is increasing. Assume that $\Gamma(G_n)\rightarrow \infty$ (otherwise the assertion trivially holds). It implies ${\rm{b}}(G_n)\rightarrow \infty$. Hence, for each $n\geq 1$, there exists an integer $N(n)$ such that ${\rm{b}}(G_i)\geq \Gamma(G_n)$ for each $i\geq N(n)$. Now, define a function $f$ by putting for each $n$, $f({\rm{b}}(G_n)):={\rm{b}}(G_{N(n)})$. We have $\Gamma(G_n)\leq f({\rm{b}}(G_n))$ for each $n$, as desired.
\end{proof}

\noindent The following result shows that the family of tree graphs is $(\Gamma, b)$-bounded.

\begin{prop}
For any tree $T$, $\Gamma(T)\leq 2{\rm{b}}(T)+2$.
\end{prop}

\noindent \begin{proof}
Set $|V(T)|=n$, $\Gamma(T)=p$ and $m(T)=m$. It is enough to show $p \leq 2m$. Otherwise, $p\geq 2m+1$. Let $p-m= m+t+1$, for some $t\geq 0$. By the definition of $m(T)$, for each $k\geq m+1$, $n-k+1$ vertices in $T$ have degree at most $k-2$. Take $k=p-m-t$ and obtain that there exist $n-p+m+t+1$ vertices of degree at most $p-m-t-2$. From the other side, there exists a Grundy coloring of $T$ using $p$ colors. Then for each $i$, at least $i$ vertices have degree at least $p-i$. Equivalently, at most $n-i$ vertices of degree at most $p-i-1$ exist in the graph. Combining these two bounds for $i=m+t+1$, we obtain $2m+2t+2\leq p$, a contradiction.   	
	\end{proof}

\noindent In the following, we denote the path on $k$ vertices by $P_k$. For any fixed graph $H$, by $Forb(H)$ we mean the family of all graphs $G$ which does not contain $H$ as induced subgraph. $Forb(H_1, H_2)$ is defined similarly.

\begin{prop}
$Forb(P_k)$ is $(\Gamma, b)$-bounded if and only if $k\leq 5$.\label{path}
\end{prop}

\noindent \begin{proof}
Define a bipartite graph $B_t$, $t\geq 2$ as follows. Take a complete bipartite graph $K_{t,t}$ and remove the edges of a matching of size $t-1$ from the graph and call it $B_t$. It's easily seen that $\Gamma(B_t)=t+1$. It can also be shown that ${\rm{b}}(B_t)=2$. Note that $B_t$ contains $P_5$ as induced subgraph but not $P_k$ for each $k\geq 6$ and hence is $P_k$-free for each $k\geq 6$. Therefore, the family of $P_k$-free graphs is not $(\Gamma, b)$-bounded for $k\geq 6$.

\noindent Assume now that $G$ is any $P_5$-free graph. A result of Kierstead et al. \cite{KPT} asserts that the family of $P_5$-free graphs is $(\chi_{\sf FF}, \omega)$-bounded. It follows that the very family is $(\Gamma, b)$-bounded.
\end{proof}

\noindent We say a graph $G$ is ${\rm b}$-monotone if for each induced subgraph $H$ of $G$ we have ${\rm{b}}(H)\leq {\rm{b}}(G)$. The family of non ${\rm b}$-monotone graphs is not $(\Gamma, b)$-bounded. For this purpose it's enough to consider the graphs $B_t$, $t\geq 2$, introduced in the proof of Proposition \ref{path}. Recall that $\Gamma(B_t)=t+1$ but ${\rm{b}}(B_t)=2$, for each $t\geq 2$. Also $B_t$ is not ${\rm b}$-monotone for each $t\geq 4$, because by removing the two vertices of degree $t$ in $B_t$ we obtain a subgraph with ${\rm b}$-chromatic number $t-1$. We make the following conjecture.

\begin{conj}
There exists a function $f(x)$ such that if $G$ is any ${\rm b}$-monotone graph then $\Gamma(G)\leq f({\rm{b}}(G))$.
\end{conj}

\noindent The next proposition proves that the conjecture is valid for all $K_{t,t}$-free graphs, for any fixed integer $t$. We need to define a tree $R_k$ of radius two. Take a vertex $v$ of degree $k-1$ as the root of $R_k$. Let $v_1, \ldots, v_{k-1}$ be the children of $v$. For each $i$, attach $k-2$ vertices of degree one to $v_i$. These vertices are all distinct so that $R_k$ contains $(k-1)(k-2)$ leaves. It is easily seen that $R_k$ admits a ${\rm b}$-coloring using $k$ colors, where $v_i$ is color-dominating vertex of color $i$.

\begin{prop}
Let $t\geq 2$ be any fixed integer and $\{G_n\}_{n\geq 1}$ be any sequence of $K_{t,t}$-free ${\rm b}$-monotone graphs. Then $\Gamma(G_n)\rightarrow \infty$ implies ${\rm{b}}(G_n)\rightarrow \infty$.
\end{prop}

\noindent \begin{proof}
Assume on the contrary that $\{G_n\}_{n\geq 1}$ is a sequence of $K_{t,t}$-free ${\rm b}$-monotone graphs with $\Gamma(G_n)\rightarrow \infty$ but for some integer $p$ and any $n$, ${\rm{b}}(G_n)\leq p$. Since $G_n$ is ${\rm b}$-monotone then $R_{p+1}$ is not an induced subgraph of $G_n$. Hence $G_n$ is $(K_{t,t}, R_{p+1})$-free for each $n$. A result of Kierstead and Penrice \cite{KP} asserts that $\{G_n\}_{n\geq 1}$ is $(\chi_{\sf FF}, \omega)$-bounded. Hence for some function $f(x)$, $\Gamma(G_n)\leq f(\omega(G_n))\leq f({\rm{b}}(G_n))\leq f(p)$, a contradiction.
\end{proof}

\end{document}